\documentstyle[12pt, leqno]{ article}
\begin{document}

\def\BC{\hbox{\rm\ l\hskip -2.2truemm C}}
\def\C{\hbox{\rm\ l\hskip -1.5truemm C}}
\def\P{\hbox{\rm\ l\hskip -0.7truemm P}}
\def\N{\hbox{\rm\bf N}}
\def\Z{\hbox{\rm\bf Z}}
\def\mult{{\rm mult}}
\def\ord{{\rm ord}}
\author{ Nguyen Van Chau\thanks {  Supported in part by the National Basic Program on Natural Science, Vietnam.} \\ {\small Hanoi Institute of Mathematics, P. O.  Box 631,  Boho 10000,  Hanoi,  Vietnam.} \\  {\small E-mail:{\it nvchau@thevinh.ncst.ac.vn}} }

\title{ Two remarks on non-zero constant Jacobian polynomial map of $\BC^2$\footnote {revised version of that send to Ann. Pol. Math. by post at April 9, 2002}
}

\date {}

\maketitle

\begin{abstract}
We present some estimations on geometry of the exceptional value sets of non-zero constant Jacobian polynomial maps of $\C^2$ and it's components.

2000{ \it Mathematical Subject Classification:} 14 H07,  14R15.

\end{abstract}

\medskip

{\bf 1. Introduction. } Recall that the {\it exceptional  value set } $E_h$ of  a polynomial mapping $h:\C^m\longrightarrow \C^n$ is the smallest subset $E_h\subset \C^n$ such that the restriction $h:h^{-1}(\C^n\setminus E_h)\longrightarrow \C^n\setminus E_h$ gives a locally trivial smooth fibration. The mysterious Jacobian conjecture (See [BMW] and [E]), which posed first by Keller in 1939 and still open even for two dimensional case,  asserts that a polynomial map  of $\C^n$ with non-zero constant Jacobian must be a polynomial bijection. Consider a polynomial map $f=(P,Q):\C^2_{(x,y)} \longrightarrow \C^2_{(u,v)}$ and denote $J(P,Q):=P_xQ_y-P_yQ_x$. It is well-known that if $J(P,Q)\equiv const. \neq 0$ but $f$ is not bijective, then  the sets  $E_P$ and $E_Q$ must be non empty finite sets and  $E_f$ is a curve so that each of it's irreducible  components is a polynomially parameterized curve ( the image of a polynomial map from $\C$ into $\C^2$) (See, for example in [J1, J2]).

In this note we present the followings.

\medskip

\noindent {\bf Theorem 1. }{\it Suppose $f=(P,Q)$ is a  polynomial map with non-zero constant Jacobian. Then,  a value  $\ u_0\in E_P$ if and only if the line  $\  u=u_0$ tangents to an irreducible local branch of $E_f$. }

\medskip
This theorem  leads to the following fact that may  be used  to consider the non-existence of non-zero constant Jacobian polynomial maps with exceptional value curve of given types.

\medskip
\noindent{\bf Theorem 2. }{\it The exceptional value set of a polynomial map of $\C^2$ with non-zero constant Jacobian can not be  isomorphic to a curve composed of the images of some polynomial maps of the form $t\mapsto (t^k,q(t))$, $k\in \N$, $q(t)\in \C[t]$.}

\medskip

A special property of the curve with irreducible components given by parameterizations in theorem 2 is that it's local irreducible branches may tangent  to  the only line $u=0$.  Simply connected curves are simple examples of such curves.  This can be easy deduced from  Lin-Zaidenberg's theorem on simply connected curves (Theorem B in [ZL]).

\medskip
\noindent{\bf Corollary 3. }{\it The exceptional value set of a polynomial map of $\C^2$ with non-zero constant Jacobian can not be  a simply  connected curve.}

\medskip

Proofs of Theorem 1 and Theorem 2  will be presented in \$ 3 and \$4.

\medskip

{\bf 2. Preliminaries. } In this section, we present some elementary facts which will be used in the proof of Theorem 1.

\medskip

i) We will work with finite factional series with parameter $\xi$ of the form
$$\varphi (x,\xi)=\sum_{k=1}^{K-1} a_kx^{n_k\over m}+\xi x^{n_K\over m}, 0\neq a_k\in\C, n_k \in \Z,\  m\in \N \eqno(1)$$
with $ n_1>n_2>\dots >n_K $ and $ \gcd(\{n_k: k \leq K\})=1$.

Let  $h(x,y)$ be a non-constant primitive polynomial monic in $y$, $\deg_y h=\deg h$. A fractional series $\varphi(x,\xi)$ in (1)  is called {\it a Newton-Puiseux type } of $h$ if
$$h(x,\varphi (x,\xi))=h_\varphi(\xi)+\mbox{lower terms in } x, \ h_\varphi\in \C[\xi], \  \deg h_\varphi >0.$$ Denote by $\mult(\varphi):=m$ and $\ord (\varphi):=n_K$.

The Newton-Puiseux types of $h(x, y)$ can be constructed from  Newton-Puiseux expansions  at infinity of the curve $h(x, y)=0$.  In fact,  if  $y=y(x)$ is such a Newton-Puiseux expansion at infinity,  then there is a unique Newton-Puiseux type   $\varphi$ of $h$ and a number $\xi_0\in \C$ such that $y(x)=\varphi(x, \xi_0)+\mbox{  {\rm lower terms in }} x $.

For a Newton-Puiseux type $\varphi$ of $h$ the rational map  $\Phi:\C\times \C\longrightarrow \C\P^2$ defined by
$$\Phi(t, \xi):=(t^{-\mult(\varphi)}, \varphi(t^{-\mult(\varphi)}, \xi)).\eqno (2)$$ determines a unbranched $i_\varphi$ sheeted covering from $\C^*\times \C$ onto $\Phi (\C^*\times \C)\subset \C^2$, where $i_\varphi:=\mult(\varphi)/\gcd\{k: a_k\neq 0,\ k<K\}$. Hence, one can use the polynomial $H_\varphi(t,\xi):=h\circ \Phi (t,\xi) $ as a kind of ``extension" of $h(x,y)$.

\medskip

\noindent {\bf Lemma 4.} ( See Theorem 2.4,  [C]).  {\it Suppose $h(x, y)$ is a primitive polynomial and monic in $y$.  Then,   $c\in E_h$  if and only if $c$ is either a critical value of $h$ or a critical value of  $h_\varphi(\xi) $ for any Newton-Puiseux type $\varphi $ of $h$. }

\medskip
In fact, by Newton's theorem the polynomial $h(x,y)-c$  can be factorized as
$$h(x,y)-c=C\prod_{u}(y-u(x)),\eqno (3)$$
where the product runs over all Newton-Puiseux expansion at infinity of the curve $h=c$. Substituting $y=\varphi (x,\xi)$ into this representation, one can see that a number $\alpha \in \C$ is a zero point of $h_\varphi(\xi)-c=0$ with multiplicity $n$ if and only if the level $h=c$ has exactly $n$ Newton-Puiseux expansions at infinity of the form $\varphi(x,\alpha)+\mbox{lower term in } x$. Thus, Puiseux data at infinity of
the curves $h=c$ must be changed  when $c$ is a critical values of $h_\varphi$.

\medskip

ii) Consider a polynomial map $f=(P,Q):\C^2 \longrightarrow \C^2$. By definition  the  exceptional  value set $E_f$  is composed by the critical value set of $f$ and the so-called {\it non-proper value set} $A_f$ of $f$ - the set of all value $a\in \C^2$ such that there exists a sequence $\C^2\ni \ p_i \rightarrow \infty$ with $ f(p_i) \rightarrow a.$ Following [J1], the non-proper value set $A_f$, if non-empty, is a curve composed by the images of some polynomial maps from $\C$ into $\C^2$.

A series $\varphi (x,\xi)$ in (1) is called  {\it dicritical series } of $f$ if
$$f(x, \varphi (x, \xi))=f_\varphi (\xi) +{\mbox {\rm lower terms in }} x, \  \deg f_\varphi >0.$$

\medskip
\noindent {\bf Lemma 5.}

$$ A_f = \bigcup_{\varphi  \ is\ a\ dicritical\ series\ of\  f} f_{\varphi} (\C). $$

\medskip

\noindent {\it Proof: } Let $\varphi$ be a dicritical series of $f$ and $\Phi(t,\xi)$ be as in (2). The map $\Phi$ sends $\C^*\times\C$ to $\C^2$  and the line $\{ 0\}\times \C $ to the line at infinity of $\C\P^2$. Then, the polynomial map $F_\varphi (t, \xi):=f\circ \Phi(t, \xi)$  maps the line $\{ 0\}\times \C$ into $A_f \subset\C^2. $
Therefore,  $f_\varphi(\C)$ is an irreducible component of $A_f $, since  $\deg f_\varphi >0$.

Conversely,  assume that $\ell$ is an irreducible component of $A_f. $ By definition we can choose a smooth point $(u_0,v_0)$ of $\ell$ and  an irreducible branch at infinity $\gamma$ of the curve $P=u_0$ (or the curve $Q=v_0$) such that the image $f(\gamma)$ is a branch curve  intersecting transversally  $\ell$ at $(u_0, v_0). $ Let $\varphi (x, \xi)$ be the Newton-Puiseux type of $P$ constructed corresponding to a Newton-Puiseux expansion at infinity of $\gamma$. Then, by definitions we can verify that $\varphi$ is a dicritical series of $f$ and $f_\varphi (\C)=\ell$.
Q.E.D

\medskip

{\bf 3. Proof of Theorem 1.} We consider a  polynomial map $f:\C^2 \longrightarrow \C^2$ with non-zero constant Jacobian. Fix a suitable affine coordinate $(x,y)$ so that  $P$ and $Q$  is monic in $y.$  For  a series $\varphi$ in (1) represent
$$
\begin{array}{c}
P(x, \varphi (x, \xi))=p_\varphi(\xi)x^\frac{ a_\varphi}{ \mult(\varphi)}+\mbox {\rm lower terms in } x \\
Q(x, \varphi (x, \xi))=q_\varphi (\xi)x^\frac{ b_\varphi}{ \mult(\varphi) }+\mbox {\rm lower terms in } x,
\end{array}
$$

\medskip

\noindent {\bf Lemma 6.} {\it i) Let $\varphi$ be a Newton-Puiseux type of $P$. If $\varphi$ is not a dicritical series of $f$, then  $ \deg p_\varphi (\xi) =1$, $q_\varphi (\xi)\equiv {\rm const.}  \neq 0 $ and $b_\varphi >0$.

ii) A dicritical series of $f$ must be a Newton-Puiseux type both of $P$ and $Q$.}

\medskip

{\it Proof.} i) First, we will show that
$$\deg p_\varphi=1,\  q_\varphi \equiv const\neq 0.\eqno (*)$$
Taking  differentiation of $F(t^{-\mult(\varphi)},\varphi (t^{-\mult(\varphi)},\xi ))$, as $a_\varphi=0$ and $b_\varphi \neq 0$ we have
$$\mult(\varphi) J(P,Q)t^{\ord(\varphi)-\mult(\varphi)-1}
=-b_\varphi\dot p_\varphi q_\varphi  t^{-b_\varphi-1}+ \mbox {{\it higher terms in }}t.
$$
Since $J(P,Q)\equiv const.\neq 0$ and $\deg p_\varphi >0$, we get (*).

Now, assume the contrary that $b_\varphi < 0$. Then, there exists a Newton-Puiseux root at infinity $u(x)$ of the curve $Q=0$ such that $u(x)=\varphi(x,\xi_0)+\mbox{lower term in } x$. Let $\psi (x,\xi)$ be the Newton-Puiseux type of $Q$ constructed corresponding to $u(x)$. Obviously, $a_\psi >a_\varphi=0$ and hence $\psi$ is not a dicritical series of $f$. Furthermore, $\varphi (x,\xi)=\psi(x,\alpha)+\mbox{ lower terms in } x$ for a zero point $\alpha$ of $p_\psi(\xi)$. This is impossible, since $p_\psi(\xi)\equiv const. \neq 0$ by applying (*) to the Newton-Puiseux type $\psi$ of $Q$. Thus, we get $b_\varphi >0$.

ii) This is obtained from (i) and definitions.
 Q.E.D

\medskip

{\it Proof of Theorem 1.} If $E_f=\emptyset$, then $f$ is bijective and $E_P=E_Q=\emptyset$. Hence, we need consider only the situation when $E_f \neq \emptyset$. In this situation $E_f=A_f$, since $f$ has not singularity.

First, suppose the line $u=u_0$ tangents to a local irreducible branch curve of an irreducible component $\ell$ of $A_f$. By Lemma 5, there is a dicritical series $\varphi$ of $f$ such that $\ell$ is the image of $f_\varphi :=(p_\varphi,q_\varphi)$.
By Lemma 6 (ii) the series $\varphi$ is a Newton-Puiseux type both of $P$ and $Q$. Since the line $u=u_0$ tangents to $\ell$, $u_0$ must be a critical value of $p_\varphi$.  Hence, by Lemma 4 the number $u_0$ is an exceptional value of $P$.

Conversely, let $u_0\in \C$ and denote $L:=\{u=u_0\}$. Assume that $L$ intersects transversally each of  local irreducible branches of $A_f$ located at points in $L\cap A_f$. We have to show that $u_0$ can not be an exceptional value of $P$. In view of Lemma 4 and Lemma 6 we need only to verify that  $u_0$ is a regular value of $p_\varphi(\xi)$ for every dicritical series $\varphi$  of $f$.

Let $\varphi$ be a given  dicritical series of $f$ and $\ell:=f_\varphi(\C)$, which  is an irreducible component of $A_f$ by  Lemma 5. Let $(u_0,v_0)\in \ell$ and $\xi_0\in \C$ with $f_\varphi(\xi_0)=(u_0,v_0)$. We have to show that $\dot p_\varphi (\xi_0)\neq 0$.

Consider the polynomial map $F:\C\times \C \longrightarrow \C^2$
$$F(t,\xi):=f\circ(t^{-\mult(\varphi)},\varphi(t^{-\mult(\varphi)},\xi))=(p_\varphi,q_\varphi)(\xi)+\mbox{ higher terms in  } t.$$
For this map $F(\{0\}\times \C)=\ell$ and
$$\det DF(t,\xi)=-\mult(\varphi) J(P,Q)t^{\ord(\varphi)-\mult(\varphi)-1}.$$Since $J(P,Q)\equiv const.\neq 0$, $F$ has singularity only on the line $t=0$.

Let $\gamma:=f_\varphi (\{ \xi:\vert \xi-\xi_0\vert <\epsilon\}$
for an enough small $\epsilon >0$. As assumed the line $L$
intersects transversally $\gamma$ at $(u_0,v_0)$. So, we can
choose an enough small neighborhood $U$ of $(u_0,v_0)$ so that
$\gamma:=\ell\cap U$ is a smooth branch curve parameterized  by
$v=v_0+h(u-u_0)$ for a holomorphic function $h$, $h(0)=0$. Define
new coordinates  $(\bar u, \bar v)=(u-u_0,v-v_0-h(u-u_0))$ in $U$
and $(\bar t,\bar\xi)=(t,\xi-\xi_0)$ in an enough small
neighborhood $V$ of $(0,\xi_0)$. Let $\bar F=(\bar F_1,\bar F_2)$
is the representation of $F$ in these coordinates. Then,
$$\bar F_1(\bar t,\bar \xi)= p_\varphi(\bar \xi)-u_0+\mbox{ higher terms in } \bar t,\eqno (4)$$
$\bar F(0,0)= (0,0)$, $\bar F(\{\bar t=0\})\subset \gamma=\{\bar v=0\}$ and $\det D \bar F(\bar t,\bar \xi)\neq 0$ for $\bar t \neq 0$. Then, by examining Newton diagrams of  $\bar F_1$, $\bar F_2$ and $\det D\bar F$ we can verify that
 $$\bar F(\bar t,\bar \xi)=(\bar \xi u_1(\bar t,\bar \xi)+\bar t u_2(\bar t, \bar \xi), \bar t^k v_1(\bar t,\bar \xi)),\eqno (5)$$
where $u_1$, $u_2$ and $v_1$ are holomorphic functions define in
$V$, $u_1(0,0)\neq 0$ and $v_1(0,0)\neq 0$ (See, for example [O,
Lemma 4.1]). From (4) and (5) it follows that $\dot p_\varphi
(\xi_0) \neq 0$. Q.E.D

\medskip

\noindent{\it Remark.} From Lemma 6 (ii) and Lemma 5 one can easy see that  the exceptional value set $E_f$ of  nonsingular polynomial map  $f$ can not contains an irreducible component isomorphic to a line.

\medskip
{\bf 4. Proof of Theorem 2. } Let $f=(P,Q)$ be representation of $f$ in the coordinate in which $E_f$ consists with the images of some polynomial maps of the form $t\mapsto (t^k,q(t))$, $k\in \N$. By applying Theorem 1  we have that  $E_{ P}\subset\{0\}$. As $ f$ has not singularity, $ P$ is a nonsingular primitive polynomial. Then, from Suzuki's equality
$$ \sum_{ c\in \C }(\chi_c-\chi )= 1-\chi $$
(See in [S]) we get $\chi_0=1$.  Here,  $\chi_c$ and  $\chi$ indicate  the Euler-Poincare characteristic of the fiber $ P=c$ and the generic fiber of $ P$,  respectively. Since the curve $ P=0$ is smooth,  it has one connected component $\ell$ diffeomorphic to $\C$. This component $\ell$ must be isomorphic to $\C $ by the Abhyankar-Moh Theorem [AM] and  the restriction of $f$ on $\ell$ must be injective.  Then, as  observed by Gwozdziewicz in [G], $f$ must be bijective. This is impossible,  since $E_f \neq \emptyset$.Q.E.D

\medskip

\noindent{\bf Acknowledgments.} The author wishes to thank Prof.  V. H.  Ha  and the referee for many valuable suggestions  and useful discussions.

\medskip

\noindent {\bf References}

 \noindent [AM] S.S.  Abhyankar and T.T.   Moh,     {\it  Embeddings of the line in the
plane},     J.   Reine Angew.   Math.   276 (1975),     148-166.

\noindent [BCW] H.   Bass,     E.   Connell and D.   Wright,    {\it  The Jacobian conjecture:
reduction of degree and formal expansion of the inverse},     Bull.   Amer.   Math.
Soc.   (N.S.) 7 (1982),     287-330.

\noindent [C] Nguyen Van Chau,  {\it
Non-zero constant Jacobian polynomial maps of $ C^2. $}
J.  Ann.  Pol.  Math.  71,  No. 3,  287-310 (1999).

\noindent [E] van den Essen,  Arno,  {\it
Polynomial automorphisms and the Jacobian conjecture}.  (English.  English summary) Progress in Mathematics,  190.  Birkhäuser Verlag,  Basel,  2000.

\noindent [G] J.  Gwozdziewicz,  {\it Injectivity on one line}.  Bull.  Soc.  Sci.  Lodz 7 (1993),  59-60,  Series: Recherches sur les deformationes XV.

\noindent [J1] Z.  Jelonek, {\it The set of points at which a polynomial map is not proper.} Ann.  Pol.  Math,  58 (1993),  259-266.

\noindent [J2] Z.Jelonek,{\it
Testing sets for properness of polynomial mappings}, Math.
Ann. 315, 1999, 1-35.

\noindent [O]  S.   Yu Orevkov,     {\it  On three-sheeted polynomial mappings of $C^2$
}, Izv.   Akad.   Nauk USSR 50 (1986),      No 6,      1231-1240,      in Russian.

\noindent [S] M.   Suzuki,     {\it Proprietes topologiques des polynomes de deux variables compleces at automorphismes algebriques  de lespace $\C^2$,    } J.   Math.   Soc.   Japan 26,     2 (1974),     241-257.

\noindent [ZL] M.G.  Zaidenberg and V.Ya.  Lin,  {\it An irreducible simply connected algebraic curve in $\C^2$ is equivalent to a quasihomogeneous  curve},  Soviet math.  Dokl.  Vol 28 (1983),  No 1,
200-205.

\end{document}